\documentclass[12pt]{article}
\usepackage{amsmath,amssymb,amsthm}
\usepackage{fullpage}

\newtheorem{thm}{Theorem}

\newtheorem{lemma}{Lemma}
\newcommand{\dis}{\displaystyle}

\begin{document}

\title{Sharp Bounds for the Harmonic Numbers}
\author{Mark B. Villarino\\
Depto.\ de Matem\'atica, Universidad de Costa Rica,\\
2060 San Jos\'e, Costa Rica}
\date{\today}

\maketitle

 \begin{abstract}
We obtain best upper and lower bounds for the \textsc{Lodge-Ramanujan} and \textsc{DeTemple-Wang} approximations to the nth Harmonic Number.
 \end{abstract}

\section{Introduction}

For every natural number $n\geqslant 1$ the \textbf{\emph{Harmonic Number}}, $H_{n}$ is the $n$th partial sum of the harmonic series:\begin{equation}
\fbox{$\dis H_{n}:=1+\frac{1}{2}+\frac{1}{3}+\cdots+\frac{1}{n}.$}
\end{equation}

Although the asymptotics of $H_{n}$ were determined by \textsc{Euler}, (see \cite{K}), in his famous  formula: \begin{equation}
\fbox{$\dis H_{n}\sim \ln n +\gamma +\frac{1}{2n}-\frac{1}{12n^{2}}+\frac{1}{120n^{4}}-\left[\cdots\right],$}
\end{equation}where $\gamma=0.57721\cdots$ is \textsc{Euler}'s constant and each summand in the asymptotic expansion is of the form $\dfrac{B_{k}}{n^{k}}$, where $B_{k}$ denots the $k$th \textsc{Bernoulli} number, mathematicians have continued to offer alternate approximative formulas to \textsc{Euler}'s. We cite the following formulas, which appear in order of increasing accuracy:\begin{align}
H_{n}    &\approx \ln n +\gamma +\frac{1}{2n+\frac{1}{3}}    \\
    &\approx\ln\sqrt{n(n+1)}+\gamma +\frac{1}{6n(n+1)+\frac{6}{5}}\\
    &\approx\ln\left(n+\frac{1}{2}\right)+\gamma +\frac{1}{24\left(n+\frac{1}{2}\right)^{2}+\frac{21}{5}}.  
\end{align}The formula (3) is the \textsc{T\'oth-Mare} approximation, (see \cite{TM}), and it\textbf{\emph{underestimates}} the true value of $H_{n}$ by terms of order $\dfrac{1}{72n^{3}}$; the second, (4), is the \textsc{Lodge-Ramanujan} approximation, and it \textbf{\emph{overestimates}} the true value of $H_{n}$ by terms of order $ \dfrac{19}{3150\left[n(n+1)\right]^{3}}$, (see \cite{Vill}); and the last, (5), is the \textsc{DeTemple-Wang} approximation, and it \textbf{\emph{overestimates}} the true value of $H_{n}$ by terms of order $\dfrac{2071}{806400\left(n+\frac{1}{2}\right)^{6}},$ (see \cite{D}).

In 2003, \textsc{Chao-Ping Chen} and \textsc{Feng Qi}, (see \cite{CQ}), published a proof of the following sharp form of the \textsc{T\'oth-Mare} approximation:

\begin{thm} For any natural number $n\geqslant 1$, the following inequality is valid:
\begin{equation}
\fbox{$\dis \frac{1}{2n+\frac{1}{1-\gamma}-2}\leqslant H_{n}-\ln n -\gamma <\frac{1}{2n+\frac{1}{3}}.$}
\end{equation}The constants $\frac{1}{1-\gamma}-2=.3652721\cdots$ and $\frac{1}{3}$ are the best possible, and equality holds only for $n=1.$

\end{thm}
The first \emph{statement} of this theorem had been announced ten years earlier by the editors of the ``Problems" section of the \emph{American Mathemtical Monthyly}, Vol 99, No. 7, (Jul-Aug, 1992), p 685, as  part of a commentary on the solution of Problem 3432, but they did not publish the proof.  So, the first \emph{published} proof is apparently that of \textsc{Chen} and \textsc{Qi}.

In this paper we will prove sharp forms of the \textsc{Lodge-Ramanujan} approximation and the \textsc{DeTemple-Wang} approximation.

\begin{thm} For any natural number $n\geqslant 1$, the following inequality is valid:
\begin{equation}
\fbox{$\dis \frac{1}{6n(n+1)+\frac{6}{5}}< H_{n}-\ln\sqrt{n(n+1)}-\gamma \leqslant\frac{1}{6n(n+1)+\frac{12\gamma -11-12\ln 2}{1-\gamma-\ln\sqrt{2}}}.$}
\end{equation}The constants $\frac{12\gamma -11-12\ln 2}{1-\gamma-\ln\sqrt{2}}=1.12150934\cdots$ and $\frac{6}{5}$ are the best possible, and equality holds only for $n=1.$
\end{thm} 

\noindent and

\begin{thm} For any natural number $n\geqslant 1$, the following inequality is valid:
\begin{equation}
\fbox{$\dis\frac{1}{24\left(n+\frac{1}{2}\right)^{2}+\frac{21}{5}} \leqslant H_{n}-\ln n -\gamma <\frac{1}{24\left(n+\frac{1}{2}\right)^{2}+\frac{54\ln\frac{3}{2}+54\gamma-53}{1-\ln\frac{3}{2}-\gamma}}.$}
\end{equation}The constants $\frac{54\ln\frac{3}{2}+54\gamma-53}{1-\ln\frac{3}{2}-\gamma}=3.73929752\cdots\cdots$ and $\frac{21}{5}$ are the best possible, and equality holds only for $n=1.$
\end{thm}

All three theorems are corollaries of the following stronger theorem:
\begin{thm}For any natural number $n\geqslant 1$, define $f_{n}$, $\lambda_{n}$, and $d_{n}$ by:\begin{align}
H_{n}    &:= \ln n +\gamma +\frac{1}{2n+f_{n}}    \\
    &:=\ln\sqrt{n(n+1)}+\gamma +\frac{1}{6n(n+1)+\lambda_{n}}\\
    &:=\ln\left(n+\frac{1}{2}\right)+\gamma +\frac{1}{24\left(n+\frac{1}{2}\right)^{2}+d_{n}},  
\end{align}respectively. Then for any natural number $n\geqslant 1$ the sequence $\{f_{n}\}$  is \textbf{monotonically decreasing} while the sequences $\{\lambda_{n}\}$ and $\{d_{n}\}$ are \textbf{monotonically increasing}.
\end{thm}\textsc{Chen} and \textsc{Qi}, (see \cite{CQ}), proved that the sequence $\{f_{n}\}$ \textbf{\emph{decreases}} monotonically.  In this paper we will prove the monotonicity of the sequences $\{\lambda_{n}\}$ and $\{d_{n}\}$.
\section{Lemmas}

Our proof is based on inequalities satisfied by the \textbf{digamma} function, $\Psi(x)$:
\begin{equation}
\fbox{$\dis \Psi(x):=\frac{d}{dx}\ln\Gamma(x)\equiv\frac{\Gamma'(x)}{\Gamma(x)}\equiv -\gamma-\frac{1}{x}+x\sum_{n=1}^{\infty}\frac{1}{n(x+n)} ,$}
\end{equation}which is the generalization of $H_{n}$ to the real variable $x$ since $\Psi(x)$ and $H_{n}$ satisfiy the equation:\begin{equation}\Psi(n+1)=H_{n}-\gamma.\end{equation} 

\begin{lemma}For every $x>0$ there exist numbers $\theta_{x}$ and $\Theta_{x}$, with $0<\theta_{x}<1$ and $0<\Theta_{x}<1$, for which the following equations are true:\begin{align}
 \Psi(x+1)   &=\ln x+\frac{1}{2x}-\frac{1}{12x^{2}}+\frac{1}{120x^{4}}-\frac{1}{252x^{6}}+\frac{1}{240x^{8}}\theta_{x}, \\
 \Psi'(x+1)   &=\frac{1}{x}-\frac{1}{2x^{2}}+\frac{1}{6x^{3}}-\frac{1}{30x^{5}}+\frac{1}{42x^{7}}-\frac{1}{30x^{9}}\Theta_{x}. \\
      \end{align}
\end{lemma}
\begin{proof}

Both formulas are well-known.  See, for example, \cite{Ed}, pp 124-125.

\end{proof}

\begin{lemma}The following inequalities are true for $x>0$:\begin{multline}
\frac{1}{3x(x+1)}-\frac{1}{15x^{2}(x+1)^{2}}< 2\Psi(x+1)-\ln\{x(x+1)\}\\
<\frac{1}{3x(x+1)}-\frac{1}{15x^{2}(x+1)^{2}}+\frac{8}{315x^{3}(x+1)^{3}}  ,   \end{multline}\begin{multline}
  \frac{1}{x^{2}}-\frac{1}{x(x+1)}-\frac{1}{3x^{3}}+\frac{1}{15x^{5}}-\frac{1}{18x^{7}}<\frac{1}{x}+\frac{1}{x+1}-2\Psi'(x+1)\\  <\frac{1}{x^{2}}-\frac{1}{x(x+1)}-\frac{1}{3x^{3}}+\frac{1}{15x^{5}}. 
   \end{multline}

\end{lemma}
\begin{proof}

The inequalities (17) were proved in our paper, (see\cite{Vill}), for integers $n$ instead of the real variable $x$.  But the proofs are valid for real $x$.

For (18) we start with (15) of \textbf{Lemma 1.}  We conclude that $$\frac{1}{2x^{2}}-\frac{1}{6x^{3}}+\frac{1}{30x^{5}}-\frac{1}{36x^{7}}<\frac{1}{x}-\Psi'(x+1)<\frac{1}{2x^{2}}-\frac{1}{6x^{3}}+\frac{1}{30x^{5}}.$$Now we multiply  to all three components of the inequality by 2 and add $\dfrac{1}{x+1}-\dfrac{1}{x}$ to them.

\end{proof}

\begin{lemma}The following inequalities are true for $x>0$:\begin{multline}
\frac{1}{\left(x+\frac{1}{2}\right)}-\frac{1}{x}+\frac{1}{2x^{2}}-\frac{1}{6x^{3}}+\frac{1}{30x^{5}}-\frac{1}{42x^{7}}< \frac{1}{x+\frac{1}{2}}-\Psi'(x+1)\\
<\frac{1}{\left(x+\frac{1}{2}\right)}-\frac{1}{x}+\frac{1}{2x^{2}}-\frac{1}{6x^{3}}+\frac{1}{30x^{5}},     \end{multline}
\begin{multline}
\frac{1}{24x^{2}}-\frac{1}{24x^{3}}+\frac{23}{960x^{4}}-\frac{1}{160x^{5}}-\frac{11}{8064x^{6}}-\frac{1}{896x^{7}}< \Psi(x+1)-\ln\left(x+\frac{1}{2}\right)\\
<\frac{1}{24x^{2}}-\frac{1}{24x^{3}}+\frac{23}{960x^{4}}-\frac{1}{160x^{5}}-\frac{11}{8064x^{6}}-\frac{1}{896x^{7}}+\frac{143}{30720x^{8}} .    \end{multline}
\end{lemma}
\begin{proof} Similar to the proof of \textbf{Lemma 2.}
\end{proof}
\section{Proof for the Lodge-Ramanujan approximation}
\begin{proof}
We solve (10) for $\lambda_{n}$ and use (13) to obtain $$\lambda_{n}=\frac{1}{\Psi(n+1)-\ln\sqrt{n(n+1)}}-6n(n+1).$$Define \begin{equation}
\fbox{$\dis \Lambda_{x}:=\frac{1}{2\Psi(x+1)-\ln x(x+1)}-3x(x+1). $}
\end{equation}for all $x>0$.  Observe that $2\Lambda_{n}=\lambda_{n}.$
\\
\\
\emph{We will show that} $\Lambda_{x}'>0$ for $x>5.$  Computing the derivative we obtain$$\Lambda_{x}'=\frac{\frac{1}{x}+\frac{1}{x+1}-\Psi'(x+1)}{\{ 2\Psi(x+1)-\ln\{x(x+1)\}^{2}}-(6x+3)$$ and therefore \begin{align*}
\{ 2\Psi(x+1)-\ln\{x(x+1)\}^{2}\Lambda_{x}'&=\frac{1}{x}+\frac{1}{x+1}-\Psi'(x+1)-(6x+3)\{ 2\Psi(x+1)-\ln\{x(x+1)\}^{2}.\end{align*}  By \textbf{Lemma 2}, this is greater than\begin{align*}
&\frac{1}{x^{2}}-\frac{1}{x(x+1)}-\frac{1}{3x^{3}}+\frac{1}{15x^{5}}-\frac{1}{18x^{7}}\\
&-(6x+3)\left\{\frac{1}{3x(x+1)}-\frac{1}{15x^{2}(x+1)^{2}}+\frac{8}{315x^{3}(x+1)^{3}}\right\}^{2}\\
&=\frac{1071x^{6}+840x^{5}-17829x^{4}-49266x^{3}-502999x^{2}-22178x-3675}{66150x^{7}(x+1)^{6}}\\
&=\frac{(x-5)\left(x^{5}+\frac{295}{51}x^{4}+\frac{628}{51}x^{3}+\frac{784}{51}x^{2}+\frac{32021}{1071}x
+\frac{137927}{1071}\right)+\frac{685960}{1071}}{\frac{1051}{17}x^{7}(x+1)^{6}}
\end{align*}which is obviously\emph{ positive} for $x>5.$  

For $x=1, \ 2, \ 3, \ 4, \ 5,$ we compute directly:\begin{align*}
\label{}
   \Lambda_{1} &=.56075467\cdots   \\
   \Lambda_{2} &=.58418229\cdots   \\
   \Lambda_{3} &=.59158588\cdots   \\
   \Lambda_{4} &=.59481086\cdots   \\
   \Lambda_{5} &=.59649019\cdots   \\
\end{align*}Therefore, the sequence $\{\Lambda_{n}\}$, $n \geqslant 1$, is a strictly increasing sequence, and therefore so is the sequence $\{\lambda_{n}\}$.

Moreover, in \cite{Vill}, we proved that $$\lambda_{n}=\frac{6}{5}-\Delta_{n},$$where $0<\Delta_{n}<\dfrac{38}{175n(n+1)}$. Therefore $$\lim_{n\rightarrow\infty}\lambda_{n}=\frac{6}{5}.$$  This completes the proof.
\end{proof}

\section{Proof for the DeTemple-Wang Approximation}

\begin{proof}
Following the idea in the proof of the  \textsc{Lodge-Ramanujan} approximation we solve (11) for $d_{n}$ and define the corresponding real-variable version.  Let
\begin{equation}
\fbox{$\dis d_{x}:=\frac{1}{\Psi(x+1)-\ln\left(x+\frac{1}{2}\right)}-24\left(x+\frac{1}{2}\right)^{2}$}
\end{equation}We compute the derivative, ask\emph{ when it is \textbf{positive}}, clear the denominator and observe that we have to solve the inequality:$$\left\{\frac{1}{x+\frac{1}{2}}-\Psi'(x+1)\right\}-48\left(x+\frac{1}{2}\right)\left\{\Psi(x+1)-\ln\left(x+\frac{1}{2}\right)\right\}^{2}>0.$$By \textbf{Lemma 3}, the left hand side of this inequality is\begin{align*}
 &>\frac{1}{\left(x+\frac{1}{2}\right)}-\frac{1}{x}+\frac{1}{2x^{2}}-\frac{1}{6x^{3}}+\frac{1}{30x^{5}}-\frac{1}{42x^{7}}-48\left(x+\frac{1}{2}\right)\\
 &\left(\frac{1}{24x^{2}}-\frac{1}{24x^{3}}+\frac{23}{960x^{4}}-\frac{1}{160x^{5}}-\frac{11}{8064x^{6}}-\frac{1}{896x^{7}}+\frac{143}{30720x^{8}}\right)^{2}   \\
 \end{align*}for all $x>0.$ This last quantity is equal to
\begin{align*}&(-9018009-31747716 x-14007876 x^2+59313792 x^3+
11454272 x^4-129239296 x^5+119566592 x^6\\
&+65630208 x^7-701008896 x^8-534417408 x^9+
178139136 x^{10})/(17340825600 x^{16} (1+2 x))\end{align*}

\noindent The denominator, $$17340825600 x^{16} (1+2 x),$$ is evidently \emph{positive} for $x>0$  and the  \emph{numerator} can be written in the form $$p(x)(x-4)+r$$ where \begin{align*}
p(x)&=548963242092+137248747452 x+34315688832 x^2
+8564093760 x^3+2138159872 x^4\\&+566849792 x^5
+111820800 x^6+11547648 x^7+178139136 x^8+178139136 x^9
\end{align*}with remainder $r$ equal to
$$r=2195843950359.$$

Therefore, the numerator is clearly \emph{positive} for $x>4,$ and therefore, the derivative, $d_{x}\ '$, too, is \emph{postive} for $x>4.$  Finally\begin{align*}
 d_{1} &=3.73929752\cdots   \\
 d_{2} &=4.08925414\cdots   \\  
 d_{3} &=4.13081174\cdots   \\
 d_{4} &=4.15288035\cdots
       \end{align*}Therefore $\{d_{n}\}$ is an \textbf{\emph{increasing}} sequence for $n\geqslant 1.$ 

Now, if we expand the formula for $d_{n}$ into an asymptotic series in powers of $\dfrac{1}{\left(n+\frac{1}{2}\right)}$, we obtain$$d_{n}\sim \frac{21}{5}-\frac{1400}{2071\left(n+\frac{1}{2}\right)}+\cdots$$and we conclude that $$\lim_{n\rightarrow\infty}d_{n}=\frac{21}{5}.$$This completes the proof.
\end{proof}

\end{document}